\documentclass[11pt,a4paper]{amsart}

\newtheorem{prop}{Proposition}
\newtheorem{thm}[prop]{Theorem}
\newtheorem{coroll}[prop]{Corollary}
\newtheorem{lemma}[prop]{Lemma}

\theoremstyle{definition}
\newtheorem{rmk}[prop]{Remark}

\newcommand{\erre}{\mathbb{R}}

\newcommand{\E}{\mathbb{E}}

\renewcommand{\P}{\mathbb{P}}

\newcommand{\F}{\mathcal{F}}

\newcommand{\be}{\begin{equation}}
\newcommand{\ee}{\end{equation}}

\renewcommand{\epsilon}{\varepsilon}
\renewcommand{\phi}{\varphi}

\newcommand{\tr}{\mathop{\mathrm{tr}}\nolimits}
\newcommand{\cp}[2]{\langle #1,#2\rangle}

\newcommand{\ds}{\displaystyle}

\begin{document}

\title[Drift reconstruction for SDEs]{On the reconstruction of the
  drift of a diffusion from transition probabilities which are partially
observed in space}

\author{S. Albeverio} 
\address{Sergio Albeverio, Institut f\"ur Angewandte Mathematik,
  Universit\"at Bonn, Wegelerstr. 6, D-53115 Bonn (Germany); IZKS; BiBoS;
  Dipartimento di Matematica, Universit\`a di Trento; CERFIM (Locarno).}
\email{albeverio@uni-bonn.de}

\author{C. Marinelli} \address{Carlo Marinelli, Institut f\"ur
  Angewandte Mathematik, Universit\"at Bonn, Wegelerstr. 6, D-53115
  Bonn (Germany).}
\email{marinelli@wiener.iam.uni-bonn.de}

\thanks{This work was carried out while the second author was on leave
  at the Institute of Mathematics (Kyiv), supported by a CNR/NATO
  grant.}

\subjclass[2000]{Primary 62M99, 60J60; Secondary 60J35, 35R30, 44A12}

\keywords{Inverse problems, stochastic differential equations, X-ray
  transform, Schr\"odinger operators, elliptic operators.}

\date{30 October 2004}

\begin{abstract}
The problem of reconstructing the drift of a diffusion in
$\erre^d$, $d\geq 2$, from the transition probability density observed
outside a domain is considered. The solution of this problem also solves a
new inverse problem for a class of parabolic partial differential equations.
This work considerably extends
\cite{jsp} in terms of generality, both concerning assumptions on the drift
coefficient, and allowing for non-constant diffusion coefficient. Sufficient
conditions for solvability of this type of inverse problem for $d=1$ are
also given.
\end{abstract}

\maketitle

\section{Introduction}
\noindent
Let $(x_t,\P^x)$ be the weak solution of the stochastic equation
\begin{equation}
\label{eq:1}
x_t = x+ \int_0^t c(x_s)\,ds + \int_0^t \sqrt{a(x_s)}\,dw_s \\
\end{equation}
where $c$ is a measurable vector field on $\erre^d$, $d\geq 2$, $a$ is
symmetric positive-definite matrix, and $w_\cdot$ is a standard
$d$-dimensional Brownian motion. Denoting with $p(x,t,y)$ the
transition probability density of $(x_t,\P^x)$, our aim is to
reconstruct $c$ from the observation of $p(x,t,y)$ for
$(x,y)\in\Lambda^c\times\Lambda^c$, $t \in [0,T]$, where $\Lambda$ is
a bounded domain in $\erre^d$ and $T>0$. We point out that this is an
entirely different problem from the ones which are usually studied in
filtering theory (or statistics of processes where the observations are
partial with respect to time instead of space).

The problem we handle in the present paper was introduced and solved in
\cite{jsp} under the assumptions that $a(x)$ is constant,
$c=\nabla\varphi$ with $\varphi\in C^2(\erre^d)$, and under growth
conditions on $c$ such that a strong solution to (\ref{eq:1}) exists with
infinite lifetime.
The results were rediscussed in \cite{ah}, where more detailed proofs can be
found, and an approach to the case of diagonal diffusion through a random
time change is sketched. We improve all these results in several directions.
In particular, we impose only integrability assumptions on $\phi$, we
require only the existence of weak solutions instead of strong ones, and we
allow $\Lambda$ to be unbounded.  Moreover, we consider much more general
classes of SDEs with variable diffusion coefficients. We also obtain a
partial solution of the problem in the case of the drift not being a
gradient field and for one-dimensional diffusions. In their full generality,
however, these latter problems were and remain unsolved (to the best of our
knowledge).

The problem at hand admits a purely analytic interpretation. Namely,
Kolmogorov's classical work implies that $p(x,t,\cdot)$ solves the
parabolic partial differential equation (PDE)
\begin{equation}
\label{eq:kolmogorov}
\left\{
\begin{array}{l}
\ds \frac{\partial p}{\partial t} = L_{c,a}^*p, \quad t>0 \\[8pt]
p(x,0,\cdot) = \delta_x(\cdot),
\end{array}
\right.
\end{equation}
where $L_{c,a}^*$ is the formal adjoint of the operator
\begin{eqnarray*}
L_{c,a}u(x) &=& \frac12 \tr\big(a(x)D^2u(x)\big) + \cp{c(x)}{Du(x)} \\
&=& \frac12 a^{ij}(x)\frac{\partial^2}{\partial x_i\partial x_j}u(x)
+ c^i(x)\frac{\partial}{\partial x_i}u(x)
\end{eqnarray*}
(here and throughout the paper we adopt the convention of summation over
repeated indices). Therefore our problem can be reformulated as an inverse
problem for PDEs: from the knowledge of the solution of
(\ref{eq:kolmogorov}) in $\Lambda^c$ alone, determine the coefficient $c$ of
the first-order term (as a function on $\erre^d$). To the best of our
knowledge the problem has not been addressed in the literature on inverse
problems for PDEs (see e.g.
\cite{belov}, \cite{isakov}), and our solution can also be seen as a
probabilistic solution to this analytic problem.

The paper is organized as follows: we collect basic assumptions,
definitions, and some known results in section \ref{sec:prelim}. In
section 3 we derive a representation formula for the transition
probabilities of diffusions whose generators satisfy certain
conditions (this class, in particular, contains distorted Brownian
motion -- see e.g. \cite{AKR-proc}). We also prove some consequences
of this representation formula that allow us to reconstruct from the
transition densities of a diffusion a function of its drift
coefficient (in section 4), and eventually its drift (in section 5).
Section 6 deals with extensions such as unbounded $\Lambda$ and
one-dimensional diffusions.

\section{Preliminaries}\label{sec:prelim}
\noindent
Unless otherwise stated, we shall work under the following standing
assumptions:
\begin{itemize}
\item[(i)] the operator $L_{c,a}$ is uniformly elliptic, i.e.
$$
\cp{a(x)\xi}{\xi} \geq \delta |\xi|^2 \quad \forall x,\,\xi \in\erre^d,
$$
for some constant $\delta>0$;
\item[(ii)] equation (\ref{eq:1}) admits a unique weak solution in $\erre^d$;
\item[(iii)] the transition probability measures of all considered diffusions
  admit a continuous density with respect to Lebesgue measure.
\end{itemize}
\begin{rmk}
  By a classical result of Zvonkin and Krylov \cite{ZK}, a weak
  solution of equation (\ref{eq:1}) exists and is unique if the
  coefficients are bounded and $a$ is continuous. A detailed study of
  conditions implying that a (generalized) diffusion has continuous
  probability densities can be found in \cite{BKR} (see also
  \cite{portenko} and \cite{BDPR}).
\end{rmk}

We shall make use of some Banach spaces. We denote by
$L^p:=L^p(\erre^d;\erre)$, for $p\in[1,+\infty[$, the space of
functions $f:\erre^d\to\erre$ such that
$$
|f|^p_{L^p} := \int_{\erre^d} |f(x)|^p\,dx < \infty.
$$
For $s\in\erre$,
let $H^{p,s}=(1-\Delta)^{-s/2}L^p$ be the usual space of Bessel
potentials on $\erre^d$. The norm in $H^{p,2}$ can be taken to be
equivalent to the one of $W^{p,2}$, the usual Sobolev space of
function with (generalized) derivatives up to order $2$ in $L^p$.
$H^{p,s}_{loc}$ is the space of functions $f$ such that $f\zeta \in
H^{p,s}$ for all $\zeta\in C^\infty_0(\erre^d)$. The following
Sobolev embedding theorem holds: $f \in
H^{p,s}_{loc}(\erre^d)$ implies that $f \in C^\alpha(\erre^d)$,
$\alpha:=s-d/p$.
In particular $f \in H^{d,2}_{loc}$ implies $f \in C^1(\erre^d)$.

Let us briefly recall some results from pinned diffusions connected to
representation of transition probability densities of diffusions (see
\cite{LZ} and \cite{QZ} for more details).
Let $(x_t,\P^x)$ be a diffusion process on $\erre^d$, endowed with its
natural filtration $(\F_t)_{t\geq 0}$.
Define a probability measure $\P^{x,y}_T$ on $\sigma(\F_t:\, t<T)$ by
$$
\left. \frac{d\P^{x,y}_T}{d\P^x} \right|_{\F_t} =
\frac{p(x_t,T-t,y)}{p(x,T,y)} \quad \forall t<T.
$$
The diffusion process $((x_t)_{t\leq T},\P_T^{x,y})$ is called a
pinned diffusion, or a diffusion conditioned on $x_0=x$ and $x_T=y$.
\par\noindent
The following Cameron-Martin formula for pinned diffusions was proved
in \cite{QZ}.
\begin{lemma}
\label{lem:cond-diff}
Let $(x_t,\P^x)$ be a diffusion process with generator
$$
L_{b,a}f = \frac12 \tr(a(x)D^2f(x)) + \cp{b(x)}{Df(x)}.
$$
Let $c:\erre^d\to\erre$ be measurable and such that
$$
\E^x\Big[
\exp\Big( \frac12\int_0^T |a^{-1/2}(x_s)c(x_s)|^2 \,ds \Big)
\Big] < \infty.
$$
Let $p_b$ and $p_{b+c}$ be the transition density functions
corresponding to the diffusions with generators $L_{b,a}$ and $L_{b+c,a}$,
respectively. Suppose that
$$
y \mapsto \E_T^{x,y}\bigg[
\exp\Big(\int_0^T \cp{a^{-1/2}(x_t)c(x_t)}{dw_t} -
\frac12 \int_0^T |a^{-1/2}(x_t)c(x_t)|^2\,dt\Big)
\bigg]
$$
is continuous. Then one has
\begin{equation}
\label{eq:CM}
\frac{p_{b+c}(x,T,y)}{p_b(x,T,y)} = \,
\E_T^{x,y} \bigg[
\exp\Big(\int_0^T \cp{a^{-1/2}(x_t)c(x_t)}{dw_t} -
\frac12 \int_0^T |a^{-1/2}(x_t)c(x_t)|^2\,dt\Big)
\bigg].
\end{equation}
\end{lemma}

\section{Transition densities for a class of diffusions}
\noindent
From now on we assume that the hypotheses of Lemma \ref{lem:cond-diff}
are satisfied and that there exists a function $\psi:\erre^d\to\erre$,
$\psi\in H^{d,2}$, such that $\nabla\psi(x) = a^{-1}(x)c(x)$. This
assumption, even though it looks quite restrictive, contains many
important examples. Let us consider, for instance, the case of unit
diffusion coefficient, i.e. Markov processes with generator
$Lu=\frac12\Delta u + \cp{c}{\nabla u}$. Assume that $c$ is regular
and $L$ admits an infinitesimally invariant measure
$\nu(dx)=\rho(x)\,dx$, i.e. that
$$
\int_{\erre^d} Lf\,\nu(dx) = 0 \quad \forall f\in C^\infty_0(\erre^d)
$$
(see e.g. \cite{ABR} for more details on these notions). It is well
known that the diffusion is reversible if and only if $L$ is symmetric
in $L^2(\erre^d,\nu)$.  One can also prove that $L$ is symmetric in
$L^2(\erre^d,\nu)$ if and only if $2c=\nabla(\log\rho)$. Therefore,
our problem of drift reconstruction can be solved for a large class of
reversible diffusions. However, the general case of $c$ not being a
gradient field is unfortunately not within the reach of our method,
except for some special situations discussed in the last section.

In this section we shall specialize Lemma \ref{lem:cond-diff} to the class of
diffusions just introduced, and deduce some important corollaries.
\begin{prop}
One has
$$
p_{b+c}(x,t,y) = p_b(x,t,y) \exp\big(\psi(y)-\psi(x)\big) \,
\E_t^{x,y}\Big[ \exp\Big(-\int_0^t V(x_s)\,ds\Big) \Big],
$$
where $V(x):=L_{b,a}\psi(x) + \frac12|a^{1/2}(x)\nabla\psi(x)|^2$.
\end{prop}
\begin{proof}
By an application of It\^o's lemma (as in \cite{K69}) we get
\begin{eqnarray*}
\psi(x_t)-\psi(x_0) &=& 
\int_0^t \cp{\nabla\psi(x_s)}{a^{1/2}(x_s)\,dw_s}
+ \int_0^t L_{b,a}\psi(x_s)\,ds \\
&=& \int_0^t \cp{a^{-1/2}(x_s)c(x_s)}{dw_s}
+ \int_0^t L_{b,a}\psi(x_s)\,ds.
\end{eqnarray*}
Therefore
\begin{eqnarray*}
&& \int_0^t \cp{a^{-1/2}(x_s)c(x_s)}{dw_s} -
\frac12 \int_0^t |a^{-1/2}(x_s)c(x_s)|^2\,ds \\
&& = \psi(x_t)-\psi(x_0) -\int_0^t L_{b,a}\psi(x_s)\,ds - 
\frac12\int_0^t|a^{-1/2}(x_s)c(x_s)|^2\,ds \\
&& = \psi(x_t)-\psi(x_0) -\int_0^t L_{b,a}\psi(x_s)\,ds - 
\frac12\int_0^t|a^{1/2}(x_s)\nabla\psi(x_s)|^2\,ds \\
&& = \psi(x_t)-\psi(x_0) -\int_0^t V(x_s)\,ds, 
\end{eqnarray*}
By a simple rewriting of (\ref{eq:CM}) we get the desired result.
\end{proof}

If $b=0$, $a=I$, i.e. $x_t$ is Brownian motion, we recover a formula
already obtained in \cite{jsp}, although under more regularity assumptions:
\begin{coroll} 
  Assume that $c$ satisfies the hypotheses of Lemma
  \ref{lem:cond-diff} with $b=0$ and $a=I$.  Setting $V(x) =
  \frac12(|\nabla\psi(x)|^2 + \Delta\psi(x))$, one has
  $$
  p_c(x,t,y) = p_0(x,t,y) \exp\big(\psi(y)-\psi(x)\big) \,
  \E_t^{x,y}\Big[ \exp\Big(-\int_0^t V(w_s)\,ds\Big) \Big],
  $$
  where $p_0(x,t,y)$ is the transition probability density of
  Brownian bridge.
\end{coroll}
Given a domain $\Lambda \subset \erre^d$ and $x,\, y \in
\partial\Lambda$, we shall denote by $(\gamma_x^{xy})_{s\in[0,t]}$
the (straight) line joining $x$ with $y$ in time $t$, i.e. the function
$$
\gamma_s^{xy} = x+(y-x)\frac{s}{t}, \quad 0 \leq s \leq t.
$$
The following proposition will be crucial for the solution of our
problem. First, let us observe that the elementary relation
$$
\int_0^t V(\gamma_s^{xy})\,ds = t \int_0^1 V(x+(y-x)s)\,ds
$$
holds.
\begin{prop}
\label{prop:limit}
Assume that $V$ is bounded from below and satisfies
the following property:
$$
V\circ\gamma^{xy}\in L^1(\erre;\erre), \;
|\gamma_n -\gamma^{xy}|_{L^\infty(\erre;\erre^d)} \to 0 
\; \Rightarrow \;
|V\circ\gamma_n - V\circ\gamma^{xy}|_{L^1(\erre;\erre)} \to 0,
$$
where $\gamma_n$ are continuous curves in $\erre^d$ with endpoints
$x$, $y$.  Then $V\circ\gamma^{xy} \in L^1(\erre;\erre)$ for almost
all $x$, $y$ (with respect to Lebesgue measure) and one has
\begin{equation}
\label{eq:phi}
\lim_{t\to 0} \frac{p_{b+c}(x,t,y)}{p_b(x,t,y)} = e^{\psi(y)-\psi(x)}
\end{equation}
and
\begin{equation}
\label{eq:xV}
\lim_{t\to 0} \frac1t \Big( \log \frac{p_{b+c}(x,t,y)}{p_b(x,t,y)} -
(\psi(y)-\psi(x)) \Big) =
-\int_0^1 V(x+(y-x)s)\,ds.
\end{equation}
\end{prop}
\begin{proof}
It clearly follows from $\psi \in H^{d,2}$ that $V \in L^d$, hence
$V_\Lambda:=V1_\Lambda \in L^1$.  Therefore, by Fubini's theorem,
$\int_{\gamma^{xy}} V_\Lambda(s)\,ds < \infty$ for almost all $x$,
$y \in \erre^d$, as desired.
We claim that
$$
\lim_{t\to 0} 
\frac{\E_t^{x,y}\Big[\exp\Big(-\int_0^t V(x_s)\,ds\Big)\Big]}{%
\exp\Big(-\int_0^t V(\gamma^{xy}_s)\,ds\Big)} = 1.
$$
Since one has, as a consequence of $V\circ\gamma^{xy} \in L^1$ (in
this proof $L^1$ stands for $L^1(\erre;\erre)$),
$$
\lim_{t \to 0} - \int_0^t V(\gamma_s)\,ds =
- \lim_{t \to 0} t\int_0^1 V(x+(y-x)s)\,ds = 0,
$$
we just need to prove that
$$
\lim_{t\to 0} 
\E_t^{x,y}\Big[\exp\Big(-\int_0^t V(x_s)\,ds\Big)\Big] = 1.
$$
In fact, we observe that for any constant $\varepsilon>0$, the
continuity of the paths of $x_\cdot$ implies that
$$
\lim_{t \to 0} \P^{x,y}_t\Big(
|x_\cdot-\gamma^{xy}_\cdot|_{L^\infty} \geq \varepsilon \Big) 
= 0,
$$
where $L^\infty$ stands for $L^\infty([0,t];\erre^d)$.
Moreover, 
\begin{eqnarray*}
&& \E_t^{x,y}\Big[\exp\Big(
\int_0^t -V(x_s)\,ds\Big) - 1 \Big]\\
&=& \E_t^{x,y}\Big[\exp\Big( \int_0^t -V(x_s)\,ds\Big) - 1 \Big|
|x_\cdot-\gamma^{xy}_\cdot|_{L^\infty}<\varepsilon \Big] \\
&& + \; \E_t^{x,y}\Big[\exp\Big( \int_0^t -V(x_s)\,ds\Big) - 1 \Big|
|x_\cdot-\gamma^{xy}_\cdot|_{L^\infty} \geq \varepsilon \Big]
\end{eqnarray*}
Using the elementary inequality $|e^x-1| \leq 1 \vee e^x$, we can write
\begin{eqnarray*}
\Big| \exp\Big( \int_0^t -V(x_s)\,ds\Big) - 1 \Big| 
&\leq& 1 \vee \exp\Big( \int_0^t -V(x_s)\,ds\Big) \\
&\leq& 1 \vee e^{t(-\inf V)},
\end{eqnarray*}
hence
\begin{eqnarray*}
&& \lim_{t \to 0} \bigg|
\E_t^{x,y}\Big[ \exp\Big( \int_0^t -V(x_s)\,ds\Big) - 1 \Big|
|x_\cdot-\gamma^{xy}_\cdot|_{L^\infty} \geq \varepsilon \Big] \bigg| \\
&\leq& \lim_{t \to 0}
\E_t^{x,y}\bigg[\Big| \exp\Big( \int_0^t -V(x_s)\,ds\Big) - 1 \Big|\;\bigg|\,
|x_\cdot-\gamma^{xy}_\cdot|_{L^\infty} \geq \varepsilon \bigg] \\
&\leq& \lim_{t \to 0} \Big( 1 \vee e^{t(-\inf V)} \Big)\, \P^{x,y}_t\big(
|x_\cdot-\gamma^{xy}_\cdot|_{L^\infty} \geq \varepsilon \big) = 0.
\end{eqnarray*}
Similarly, using the elementary inequality $|e^x-1|\leq e^{|x|}-1$, we have
\begin{eqnarray*}
&& \lim_{t \to 0} \bigg|
\E_t^{x,y}\Big[ \exp\Big( \int_0^t -V(x_s)\,ds\Big) - 1 \Big|
|x_\cdot-\gamma^{xy}_\cdot|_{L^\infty} < \varepsilon \Big] \bigg| \\
&\leq& \lim_{t \to 0} 
\E_t^{x,y}\Big[ \exp\Big( \int_0^t |-V(x_s)+V(\gamma_s)|\,ds\Big) \Big|
|x_\cdot-\gamma^{xy}_\cdot|_{L^\infty} < \varepsilon \Big]
e^{t |V\circ\gamma^{xy}|_{L^1}} - 1 \\
&\leq& \lim_{t \to 0} e^{\delta_{\varepsilon}} 
e^{t |V\circ\gamma^{xy}|_{L^1}} - 1 = 0,
\end{eqnarray*}
where we have used the following immediate consequence of the
hypotheses: for any $\varepsilon > 0$ there exists $\delta_\varepsilon > 0$ 
such that
$$
\int_0^t |V(\gamma^{xy}_s)-V(x_s)|\,ds = |V\circ\gamma^{xy} - V\circ x|_{L^1}
< \delta_\varepsilon
$$
whenever $|x_\cdot - \gamma_\cdot^{xy}|_{L^\infty}<\varepsilon$. This
concludes the proof of the claim.  
Assertion (\ref{eq:phi}) now follows immediately from the claim just
proved and the previous proposition. Assertion (\ref{eq:xV}) follows by
$$
\frac{\E_t^{x,y}\Big[ \exp\Big( \int_0^t -V(x_s)\,ds\Big) \Big]}{%
\exp\Big( \int_0^t -V(\gamma^{xy}_s)\,ds\Big)} = 1 + o(t)
$$
for $t \to 0$, hence
\begin{eqnarray*}
&& \lim_{t\to 0} \frac1t
\Big(\log\frac{p_{b+c}(x,t,y)}{p_b(x,t,y)} - \big(\psi(y)-\psi(x)\big)\Big) \\
&=&
\lim_{t\to 0} \frac1t \log
\E_t^{x,y}\Big[ \exp\Big( \int_0^t -V(x_s)\,ds\Big) \Big] \\
&=& \lim_{t\to 0} \frac1t \log 
\Big( \exp\Big( \int_0^t -V(\gamma_s)\,ds\Big) + o(t) \Big) \\
&=& \lim_{t\to 0} \frac1t \log 
\exp\Big(t \int_0^1 -V(x+(y-x)s)\,ds\Big) \\
&=& \int_0^1 -V(x+(y-x)s)\,ds,
\end{eqnarray*}
which finishes the proof.
\end{proof}

\begin{rmk}
If $x_\cdot$ is Brownian motion, Aizenman and Simon \cite{AS}
proved that $V$ is of Kato class if and only if
$$
\lim_{t \to 0}\, \sup_x\, \E^x \Big[ \int_0^t |V(x_s)|\,ds\Big] = 0,
$$
without assuming that $V$ is lower bounded. Using Kasminskii's
lemma (see e.g. \cite{portenko}) it is then immediate to deduce that
if $V$ is of Kato class, then one also has
\begin{equation}
\label{eq:khas}
\lim_{t \to 0}\, \sup_x\, \E^x\bigg[
\exp\Big( \int_0^t |V(x_s)|\,ds \Big) \bigg] = 1.
\end{equation}
Moreover, Chung and Zhao proved that (\ref{eq:khas}) continues to hold
for any Hunt process $x_\cdot$, under the only assumption that $V$ is
of Kato class (see \cite{chung}, Proposition 3.8). We also refer to
\cite{ABM} and \cite {AM} for related results on generalized
Schr\"odinger operators and forms and associated Markov processes.
\end{rmk}

\section{Reconstruction of $V$}
\noindent
In this section we assume that the hypotheses of Lemma
\ref{lem:cond-diff} and Proposition \ref{prop:limit} are satisfied.

We shall then show that the transition probabilities of $(x_t,\P^x)$
determine the X-ray transform of $V_\Lambda$ (in the sense of, e.g.,
\cite{helgason}), which in turn yields $V_\Lambda$ by a Fourier
transform argument.
In particular, from equations (\ref{eq:phi}) and (\ref{eq:xV}) of
Proposition \ref{prop:limit} it immediately follows that
$$
F(x,y) := \int_0^1 V(x+(y-x)s)\,ds, \quad x,y \in \partial\Lambda
$$
is determined by the transition probabilities of $(x_t,\P^x)$.
Moreover, one immediately recognizes that $F$ is the X-ray
transform of $V_\Lambda$.  Let us fix some notation: we represent a
line $\gamma$ as a pair $\gamma=(\omega,z)$, where $\omega \in
\mathbb{S}^{n-1}$ is a unit vector in the direction of $\gamma$ and
$z\in\gamma\cap\omega^\perp$.  Then the line integral $\int_\gamma
f(x)\,dx$ is denoted by
$$
\hat{f}(\gamma) = \hat{f}(\omega,z) = P_\omega f(z).
$$
As argued before, since we know that $V_\Lambda \in L^1$, Fubini's
theorem implies that for each $\omega\in\mathbb{S}^{d-1}$, $P_\omega
V_{\Lambda}(z)$ is defined for almost all $z \in \omega^\perp$.
Moreover, we have, for $p \in \omega^\perp$,
$$
\tilde{V}_\Lambda (p) = 
\int_{\omega^\perp} P_\omega V_\Lambda(z) e^{i\cp{p}{z}}\,dz
$$
(a result often called slice-projection theorem, see e.g.
\cite{helgason}).  One uniquely recovers $V_\Lambda$ by taking the
inverse Fourier transform of $\tilde{V}_\Lambda(p)$. Summarizing, we
have proved the following
\begin{thm}
  The restriction of $V$ to the domain $\Lambda$ can be uniquely
  reconstructed from the transition probabilities of $(x_t,\P^x)$.
\end{thm}
\begin{rmk}
  If our only aim were to reconstruct the function $V_\Lambda$ from
  its X-ray or Radon transform, even more generality could be allowed,
  up to the situation where $V$ is a distribution. For results on
  inverting the Radon transform of a distribution, see e.g.
  \cite{gv5}, \cite{RK}, \cite{helgason}.
\end{rmk}

\section{Reconstruction of the drift}
\noindent
As in the previous section, we assume that the hypotheses of Lemma
\ref{lem:cond-diff} and Proposition \ref{prop:limit} are satisfied.

Set $u(x)=e^{\psi(x)-\psi(y)}$, where $y$ is any (fixed) point on the
boundary of $\Lambda$. Then $u$ satisfies the following Dirichlet
boundary value problem for a second-order elliptic operator:
\begin{equation}
\label{eq:bvp}
\left\{\begin{array}{ll}
\frac12 a^{ij}u_{x_ix_j} + b^iu_{x_i} = V(x)u, & x \in \Lambda\\[6pt]
u(x) = e^{\psi(x)-\psi(y)}, & x \in \partial\Lambda.
\end{array}
\right.
\end{equation}
This is easily seen as a consequence of the definition of $V$ and of the
following simple calculations:
\begin{eqnarray*}
u_{x_i} &=& u \psi_{x_i} \\
u_{x_ix_j} &=& u_{x_j}\psi_{x_i} + u\psi_{x_ix_j}
= u(\psi_{x_j}\psi_{x_i} + \psi_{x_ix_j}) \\
a^{ij}u_{x_ix_j} &=& \big(a^{ij}\psi_{x_ix_j} + |a^{1/2}\nabla\psi|^2\big) u,
\end{eqnarray*}
where the last step is justified by
$$
a^{ij}\psi_{x_j}\psi_{x_i} = \cp{a\nabla\psi}{\nabla\psi} = 
\cp{a^{1/2}\nabla\psi}{a^{1/2}\nabla\psi} = |a^{1/2}\nabla\psi|^2.
$$
If (\ref{eq:bvp}) is uniquely solvable, then we are able to recover
$\psi(x)$ for $x\in\Lambda$ uniquely. In fact we have:
\begin{prop}
\label{prop:trud}
  Suppose $a^{ij}$ are differentiable and $V \in L^\infty_+(\Lambda)$.
  Then there exists a unique solution $u \in H^{2,1}(\Lambda)$ of the
  Dirichlet problem (\ref{eq:bvp}).
\end{prop}
\begin{proof}
  Since $\psi\in C^1(\Lambda)$, as follows by Sobolev embeddings, then
  $f(x):=e^{\psi(x)-\psi(y)} \in C^1(\Lambda) \subset
  H^{2,1}(\Lambda)$. Moreover, $L_{b,a}$ is strictly elliptic and
  $b=\nabla\psi \in C(\Lambda)$, hence $b$ is bounded. We can now
  appeal to Theorem 8.9 of \cite{GT}, which yields the existence and
  uniquess of a solution to (\ref{eq:bvp}), as claimed.
\end{proof}
\begin{rmk}
  Using more general results on elliptic PDEs, one can remove the
  unpleasant assumption of $V$ being bounded, at the cost of added
  technicalities. In particular, using the existence and uniquess
  results of \cite{trudinger}, one can replace $V\in
  L^\infty(\Lambda)$ by $g\in L^{d/2}$ in the hypotheses of
  Proposition \ref{prop:trud}, where
$$
g := (a^{-1})_{ij}(\tilde{b}^i+\tilde{b}^j),\quad
\tilde{b}^j := b^j -\frac12a^{ij}_{x_i}.
$$
The details (mostly calculations) are left to the reader. The
assumption $V\geq 0$ is used to guarantee that the spectrum of the
operator $L_{b,a}-V(x)$ (considered between appropriate spaces of
integrable functions) does not contain zero. If we are willing to
accept this level of generality, sacrificing a bit of concreteness,
then we can dispense with the assumption of $V$ being positive, and
simply assume that zero is not an eigenvalue of $L_{b,a}-V(x)$. For
further details we refer to \cite{GT} and \cite{trudinger}, where a
Fredholm alternative for this type of operators is established.
\end{rmk}

If $a$ is the identity matrix, hence $c=\nabla\psi$, we can obtain stronger
results. In particular, the Dirichlet problem (\ref{eq:bvp}) reduces
to the Dirichlet problem for the time-independent Schr\"odinger
operator with Hamiltonian $-\frac12\Delta+V$:
\begin{equation}
\label{eq:bvp2}
\left\{\begin{array}{ll}
\ds \frac12\Delta u = V(x) u, & x \in \Lambda\\
u(x) = e^{\psi(x)-\psi(y)}, & x \in \partial\Lambda,
\end{array}
\right.
\end{equation}
for which there exists a rich literature. We can apply, for instance,
Theorem 4.7 of \cite{chung} (see also \cite{AS}). We shall denote by
$\tau_\Lambda$ the first exit time of Brownian motion from the domain
$\Lambda$.
\begin{thm}
  Assume that $V$ is of local Kato class, $\Lambda$ is bounded and
  regular, and
  \begin{equation}
\label{eq:gauge}
  x \mapsto \E^x\Big[
  \exp\Big(-\int_0^{\tau_\Lambda}V(w_s)\,ds\Big)\Big]
  \end{equation}
  is bounded.  Then there exists a unique solution $u\in
  C(\overline{\Lambda})$ of the boundary value problem
  (\ref{eq:bvp2}).
\end{thm}
\begin{proof}
  The conditions on $V$ and $\Lambda$ are needed in order to apply the
  above mentioned results of \cite{chung}. Moreover, since $\psi\in
  C^1$ as follows by Sobolev embeddings, and thus
  $f(x):=e^{\psi(x)-\psi(y)} \in C(\partial\Lambda)$, part (iv) of
  Theorem 4.7 in \emph{ibid.} ensures that
  $$
  u(x) = \E^x\Big[ \exp\Big(-\int_0^{\tau_\Lambda} V(w_s)\,ds\Big)
  f(w_{\tau_\Lambda})\Big]
  $$
  is the unique solution of $(-\frac12\Delta+V)u=0$ such that $u\in
  C(\overline{\Lambda})$ and $u(x)=e^{\psi(x)-\psi(y)}$ on
  $\partial\Lambda$.
\end{proof}
\begin{rmk} 
  A simple sufficient condition guaranteeing that $V$ is of Kato class
  is $V \in L^{p/2}$, $p>d/2$ (see e.g. \cite{AS}, p.~233). Therefore,
  if $\psi \in H^{p,2}$ with $p>d$, $V$ is of Kato class. On the other
  hand, we were unable to find simple sufficient conditions ensuring
  that (\ref{eq:gauge}) is bounded. Let us mention, however, that each of the
  following analytic conditions is sufficient:
\begin{itemize}
\item[(i)] $\int_0^\infty T_t1\,dt$ is bounded;
\item[(ii)] $-\frac12\Delta+V \geq 0$ (in the sense of operators), or
  equivalently:
\item[(iii)] The spectrum of $-\frac12\Delta+V$ is contained in
  $]0,+\infty[$,
\end{itemize}
where we have denoted by $T_t$ the semigroup generated by
$-\frac12\Delta+V$. For more informations see \cite{AS} and
\cite{chung}, p.~126.
\end{rmk}

Once $u$ is obtained, one recovers $\psi$ immediately, and hence $c$.
We have proved the following result
\begin{thm}
  Assume that the boundary value problem (\ref{eq:bvp}) admits a
  unique solution. Then the transition probabilities $p_{b+c}(x,t,y)$
  for $x$, $y\in\Lambda^c$, $t\in[0,T]$, $T>0$, determine $c$
  uniquely.
\end{thm}
\begin{proof}
  It is just a combination of the previous steps. In particular, one
  proceeds as follows:
\begin{enumerate}
\item Obtain the X-ray transform of $V_\Lambda$ from the transition
  probabilities $p_{b+c}(x,t,y)$;
\item Invert the X-ray transform obtaining $V_\Lambda$;
\item Solve the elliptic PDE (\ref{eq:bvp}) obtaining $\psi(x)=\log
  u(x)+\psi(y)$;
\item Obtain $c=a\nabla\psi$.
\end{enumerate}
\end{proof}

\section{Some extensions}
\subsection{$a^{-1}c$ not a gradient field}
It is clear that our approach strongly relies on the assumption that
the $a^{-1}c$ is a gradient field. When this is not the case, we can
only give a rather involved sufficient condition to reduce the problem
to a more tractable one. We shall assume for simplicity that $x_\cdot$
is a $L_{c,I}$ diffusion, without knowing a priori that $c$ is a
gradient field. We also assume that the technical assumptions
introduced so far are in place when needed.
\begin{prop}
\label{prop:nongrad}
Let $f:\erre^d\to\erre^d$ be a $C^2$ diffeomorphism mapping
$\Lambda$ into a bounded domain, and such that
\begin{equation}
\label{eq:nongrad}
\big[\nabla f(f^{-1}(x))\nabla f(f^{-1}(x))^*\big]^{-1}L_{c,I}f(f^{-1}(x)) 
= \nabla\psi(x)
\end{equation}
for some $\psi:\erre^d \to \erre$. Then the transition
probabilities of $(x_t,\P^x)$ uniquely determine the drift $c$.
\end{prop}
\begin{proof}
It\^o's formula for $f_i(x_t)$, the $i$-th component of $f$, gives
$$
f_i(x_t)-f_i(x_0) = \int_0^t \cp{\nabla f_i(x_s)}{dw_s}
+ \int_0^t L_{c,I}f_i(x_s)\,ds,
$$
that is
$$
f(x_t)-f(x_0) = \int_0^t \nabla f(x_s)\,dw_s
+ \int_0^t L_{c,I} f(x_s)\,ds,
$$
or equivalently, defining $y_\cdot=f(x_\cdot)$
\begin{equation}
\label{eq:y}
y_t = y + \int_0^t L_{c,I} f(f^{-1}(y_s))\,ds +
\int_0^t \nabla f(f^{-1}(y_s))\,dw_s.
\end{equation}
The hypotheses imply that one can reconstruct the drift and the
transition probabilities of the diffusion (\ref{eq:y}). But this is
enough to recover the transition probabilities of $x_\cdot$ as well, as the
following obvious identities show:
$$
\P\big(x_t=y|x_0=x\big) = \P\big(f^{-1}(y_t)=y|f^{-1}(y_0)=x\big)
= \P\big(y_t=f(y)|y_0=f(x)\big).
$$
It is well known that the transition probabilities of an
$L_{c,I}$-diffusion determine the drift coefficient $c$ uniquely.
\end{proof}

\subsection{Unbounded $\Lambda$}
We assumed $\Lambda$ to be bounded in order to obtain existence and
uniqueness of solutions for the Dirichlet problems (\ref{eq:bvp}) and
(\ref{eq:bvp2}). However, in some cases this assumption can be
relaxed. For instance, imposing enough boundedness of the coefficients
in (\ref{eq:bvp}), one can obtain existence and uniqueness results
without assuming that $\Lambda$ is bounded.
\begin{prop}
  Assume that $a^{ij}$ are differentiable, $e^{-|x|^2}\|a(x)\| \to 0$
  for $|x|\to\infty$, and $V$ is non-negative (in the generalized sense).
  Then there exists a unique solution of (\ref{eq:bvp}).
\end{prop}
\begin{proof}
  Let $v(x)=u(x)e^{\phi(x)}=e^{\psi(x)+\phi(x)}$,
  $\phi(x)=-|x|^2$. Since $\psi\in H^{d,2}$, hence
  $\psi\in C^1$ by Sobolev embedding, $\psi$ is bounded on $\Lambda$,
  and so is $u$. Moreover, it is immediate to show that $v \in
  H^{2,1}(\Lambda)$. Then one has
\begin{eqnarray*}
L_{b-2a\nabla\varphi,a}v =
L_{b,a}v - \cp{2a\nabla\varphi}{\nabla v},
\end{eqnarray*}
and $\nabla v = (\nabla\psi+\nabla\varphi)v$,
\begin{eqnarray*}
L_{b,a}v &=& 
\big( L_{b,a}(\psi+\varphi) + |a^{1/2}\nabla(\psi+\varphi)|^2\big)v\\
&=& \big( L_{b,a}\psi + |a^{1/2}\nabla\psi|^2 +
L_{b,a}\varphi + |a^{1/2}\nabla\varphi|^2 +
2\cp{a^{1/2}\nabla\varphi}{a^{1/2}\nabla\psi} \big)v \\
&=& \Big( \frac12\tr(a\psi_{xx}) + \cp{b+2a\nabla\varphi}{\nabla\psi} +
|a^{1/2}\nabla\psi|^2 + L_{b,a}\varphi + |a^{1/2}\nabla\varphi|^2 \Big)v,
\end{eqnarray*}
therefore
\begin{eqnarray*}
L_{b-2a\nabla\varphi,a}v &=&
L_{b,a}v - \cp{2a\nabla\varphi}{\nabla v} \\
&=& L_{b,a}v - \cp{2a\nabla\varphi}{\nabla\psi+\nabla\varphi}v \\
&=& L_{b,a}v - \cp{2a\nabla\varphi}{\nabla\psi}v -
2|a^{1/2}\nabla\varphi|^2v \\
&=& \Big(\frac12\tr\,(a\psi_{xx}) + \cp{b}{\nabla\psi} +
|a^{1/2}\nabla\psi|^2 + L_{b,a}\varphi -
|a^{1/2}\nabla\varphi|^2\Big)v \\
&=& \big(V(x)+d(x)\big)v,
\end{eqnarray*}
where $d:=L_{b,a}\varphi - |a^{1/2}\nabla\varphi|^2$.
That is, $v$ solves the equation
\begin{equation}
\label{eq:v}
L_{b-2a\nabla\varphi,a}v = \big( V(x) + d(x) \big)v,
\end{equation}
with boundary condition $v(x)=f(x)e^{\phi(x)}$ on $\partial\Lambda$.
One can apply now Theorem 8.9 of \cite{GT} to determine existence and
uniqueness of a solution of (\ref{eq:v}) in $H^{2,1}$. In fact, if
$L_{b,a}$ is strictly elliptic, so is also $L_{b-2a\nabla\varphi,a}$,
and the continuity and growth condition on $a$ imply that
$b-2a\nabla\varphi$ is bounded. Finally, there is no loss of
generality in assuming that $V\geq 0$ implies $V(x)+d(x)\geq 0$: if it
were not so, we could make $|d|$ arbitrarily small by using as
cut-off function $\phi(x)=e^{-\kappa|x|^2}$, without altering the
previous results. Then the unique solution to (\ref{eq:bvp}) is given
by $u(x)=v(x)e^{-\varphi(x)}$.
\end{proof}

In the case of unit diffusion a stronger statement can be made, as
follows by the results in Chapter 5 of \cite{chung}.
\begin{prop} 
  Let $d \geq 3$, $V \in L^1(\Lambda)$ and Kato class. If $u_f
  \not\equiv \infty$ in $\Lambda$, then the solution of
  (\ref{eq:bvp2}) is given by
$$
u_f(x) = \E^x\bigg[
\exp\Big( -\int_0^{\tau_\Lambda} V(w_s)\,ds\Big) f(w_{\tau_\Lambda}) \bigg].
$$
\end{prop}
\begin{proof}
  As before, $\psi \in H^{d,2}$ implies $\psi\in C^1$, hence that
  $\psi$ is bounded on $\partial\Lambda$, and also that
  $f(x):=e^{\psi(x)-\psi(y)}\in L_+^\infty(\partial\Lambda)$. Then by
  Theorem 5.18 and 5.19 of \cite{chung} we obtain that $u_f$ solves
  $(-\frac12\Delta+V)u=0$ and $u_f \in C_b(\overline{\Lambda})$.
  Finally, $u_f$ satisfies the boundary condition as an immediate
  consequence of its definition and continuity in the closure of
  $\Lambda$.
\end{proof}

\subsection{One-dimensional case}
Using the considerations of the previous subsections it is possible to
give a solution, at least in some special cases, to the problem posed
in \cite{jsp} of reconstructing the drift of a one-dimensional
diffusion with $a=1$.  In particular, let us assume that the
transition probability density $p_1(x,t,y)$ of the diffusion
$$
dx_1(t) = c_1(x_1(t))\,dt + dw_1(t)
$$
is known for $x$, $y\in]0,1[^c$, $t\in[0,T]$, with $T$ a fixed
positive constant. As before, our aim is to determine $c_1(x)$ for
$x\in[0,1]$.

\noindent
Define an $\erre^2$-valued diffusion as weak solution of the following
SDE:
\begin{equation}
\label{eq:coupled}
\left\{
\begin{array}{l}
\ds dx_1(t) = c_1(x_1(t))\,dt + dw_1(t)\\[8pt]
\ds dx_2(t) = c_2(x_2(t))\,dt + dw_2(t),
\end{array}
\right.
\end{equation}
where $c_2$ is a smooth, bounded function, and $w_1$, $w_2$ are
standard independent Brownian motions. 
In more compact notation, we can write
$$
dx(t) = c(x(t))\,dt + dw(t),
$$
with $c(x_1,x_2) = (c_1(x_1),c_2(x_2))$, $w=(w_1,w_2)$.
It is clear that one can recover the transition probabilities of $x_1$
from those of $x$, since $x_1$ and $x_2$ are independent.

\noindent
Trying to reconstruct the vector field $c$ (which is trivially a
gradient field: $c=\nabla\psi$, $\psi(x_1,x_2)=\int_0^{x_1}c_1(s)\,ds
+ \int_0^{x_2}c_2(s)\,ds)$, one is faced with two problems:
$\Lambda=]0,1[\times\erre$ is unbounded, and $V_\Lambda$ is not in
$L^1$. One can try, however, to consider the problem of drift
reconstruction for $y_\cdot=f(x_\cdot)$, with $f$ a $C^2$
diffeomorphism.
\begin{prop} 
  Assume that there exist $c_2:\erre\to\erre$ and a diffeomorphism
  $f:\erre^2\to\erre^2$ satisfying (\ref{eq:nongrad}) and such that
  the corresponding Dirichlet problem (\ref{eq:bvp}) admits a unique
  solution on $f(\Lambda)$. Then $c_1$ can be uniquely reconstructed
  from the transition probabilities of $x_1(\cdot)$.
\end{prop}
\begin{proof} 
  The transition probabilities of $x_1(\cdot)$ uniquely determine the
  transition probabilities of $x(\cdot)$ outside the (unbounded)
  domain $\Lambda =\, ]0,1[ \times \erre$. If $f$ exists such that
  (\ref{eq:nongrad}) holds, then the the problem of reconstructing the
  drift of $y(\cdot):=f(x(\cdot))$ is well posed, and it is equivalent
  (under the technical assumptions introduced in the previous
  sections) to the solvability of the Dirichlet problem (\ref{eq:bvp})
  on the domain $f(\Lambda)$. Therefore, assuming the latter problem admits a
  unique solution, this yields the transition probabilities of
  $y(\cdot)$, hence those of $x(\cdot)$ because $f$ is a bijection. As
  already observed, the transition probabilities of a diffusion with
  generator $L_{c,I}$ uniquely determine $c$.
\end{proof}

It is clear that the last proposition is not constructive and simply gives
sufficient conditions for the solvability of the problem of drift
reconstruction in dimension one. As in the case of higher dimensional
diffusions with $a^{-1}c$ not being a gradient field, these sufficient
conditions seem difficult to check. However, since we essentially rely on
the above described representation with the drift being a gradient field,
this seems to be the best we can achieve by our present method.

\def\cprime{$'$}
\providecommand{\bysame}{\leavevmode\hbox to3em{\hrulefill}\thinspace}
\providecommand{\MR}{\relax\ifhmode\unskip\space\fi MR }
\providecommand{\MRhref}[2]{%
  \href{http://www.ams.org/mathscinet-getitem?mr=#1}{#2}
}
\providecommand{\href}[2]{#2}

\end{document}